\documentclass[oneside,english]{amsart}
\usepackage[T1]{fontenc}
\usepackage[latin1]{inputenc}
\usepackage{graphicx}

\makeatletter

\providecommand{\LyX}{L\kern-.1667em\lower.25em\hbox{Y}\kern-.125emX\@}

 \theoremstyle{plain}    
 \newtheorem{thm}{Theorem}[section]
 \numberwithin{equation}{section} %% Comment out for sequentially-numbered
 \numberwithin{figure}{section} %% Comment out for sequentially-numbered
 \theoremstyle{plain}
 \newtheorem{cor}[thm]{Corollary} %%Delete [thm] to re-start numbering
 \theoremstyle{plain}
 \newtheorem{lem}[thm]{Lemma} %%Delete [thm] to re-start numbering
 \theoremstyle{definition}
  \newtheorem{defn} [thm] {Definition}
  \newtheorem{problem}[thm]{Problem}
  \newtheorem{conj} [thm]{Conjecture}
 \theoremstyle{plain}    
 \newtheorem{prop}[thm]{Proposition} %%Delete [thm] to re-start numbering
 \theoremstyle{remark}
 \newtheorem*{rem*}{Remark}
 \newtheorem{xmpl}[thm]{Example}
 \newtheorem*{acknowledgement*}{Acknowledgement}

\usepackage{graphicx}

\usepackage{graphicx}

\usepackage{babel}
\makeatother
\begin{document}

\title{Estimates for the minimal crossing number}

\author{Hermann Gruber}

\begin{abstract}
First, I give an elementary proof for the fact
that the minimal crossing number is additive under composition of
torus links. This result is generalized to the composition of
homogeneous braids with alternating fibered links. Then there
follow estimates for the crossing number of satellite knots.
In the last chapter, I discuss a conjecture concerning the HOMFLY
and the Kauffman polynomial.
\end{abstract}
\maketitle

\section{Introduction.}

A link is called composite iff there is a 2-sphere in $S^{3}$ which
meets the link in exactly two points and decomposes it into two sublinks
K and L, neither of which is an unknotted arc. Write K\#L if the link
is the composition of K and L. Whether the minimal crossing number
of a composite link is simply the sum of the minimal crossing numbers
of the factor links, is a very natural and very old question already
posed by Tait some hundred years ago. It is trivial to see that it
cannot be greater than the sum, but the question if equality holds
in general remains unsolved up to now. However, the question is answered
if the factor links are all from a certain family of knots, called
adequate knots\cite{Kirby}, problem 1.65. Here, another such family 
is exposed, which contains
the torus links and the alternating fibered links. Some of the arguments 
exposed
here are closely following the ideas of Murasugi's paper \cite{MuBraid},
and some
results are re-expressed and proved more shortly here. It turns out to be 
astonishingly easy to show that the crossing number of torus links is 
additive under composition. But no one has written up these things 
up to now, and it is one purpose of this paper to do this.

Another aim is to turn around the classical points of view, gaining
beautiful formulas estimating the minimal crossing number via the
geometric concepts of braid index and (canonical) genus. For example, you 
will see that, in this application here, the canonical genus can 
be a considerably more powerful concept than the Seifert genus.

\section{Preliminaries.}

Throughout this paper, it is immaterial if the knots are chiral or
not. For this cause, don't regard chiral pairs as distinct.

For a (nontrivial and non-split) link K and a regular diagram D(K)
(or simply D), let c(D) be the number of crossings in D. Write c(K)
for the minimal crossing number, and b(K) for the braid index of K.

Give the knot K an arbitrary orientation. By cutting out each crossing,
respecting the orientation, convert the diagram D into a number of
oriented closed curves in the (extended) plane, called Seifert circles of D.
Write s(D) for the number of Seifert circles of D. From these circles,
construct a spanning surface for K. The genus of this surface depends
of the diagram: $g(D)=\frac{1}{2}(c(D)-s(D)-|K|+2)$, where |K| is
the number of components of K. This can easily be seen by calculating
the Euler characteristics of the constructed surface.

The (Seifert) genus of a link K is defined as the minimal genus among
all orientable surfaces bounded by K, the canonical (or weak) genus
as \[
\widetilde{g}(K)=\min _{D=D(K)}\{g(D)\}\]
 and the free genus $g_{f}(K)$ is the minimum genus of all Seifert
surfaces for K whose complement in $S^{3}$ is a handlebody. We mention
$g(K)\leq g_{f}(K)\leq \widetilde{g}(K)$.

\section{Torus links, alternating fibered links and homogeneous braids.}

I shall note first the following self-evident lemma, which plays a
central role in our discussion:

\begin{lem}
For any link K, \[c(K)=\min _{D=D(K)}\left(2g(D)+s(D)\right)+|K|-2\].
\end{lem}
By minimizing both the number of Seifert circles and the genus independently,
we get the following estimate:

\begin{lem}\label{bound1}
For every link K, \[
c(K)\geq 2\widetilde{g}(K)+b(K)+|K|-2\geq 2g(K)+b(K)+|K|-2\]

\end{lem}
\begin{proof}
What remains to show is \[\min _{D=D(K)}\{s(D)\}=b(K)\]
See Yamada's proof \cite{Yamada}.
\end{proof}

Both of the above inequalities are sharp for torus links, so we can
give a short and elementary proof of a result of Murasugi:

\begin{prop}
\cite{MuBraid}Let K be a (p,q) torus link with $p\geq q\geq 2$.
Then \[
c(K)=2g(K)+b(K)+|K|-2=pq-p\]
\end{prop}
\begin{proof}
Consider the standard representation D of K as a closed braid with
q strands. Here, $c(D)=pq-p$ and $2g(D)=pq-p-q+|K|-2$.
g(D) is minimal over all Seifert surfaces, see  \cite{Cro}. 
As the diagram D shows, the braid index can
be at most q, and since it cannot be lower than the bridge number,
it must be equal to q.
\end{proof}

This result motivates the following definition:
\begin{defn}
Say that a link K is in the family $\mathcal{F}$ 
(or, shortly, $\mathcal{F}$-link) if
$c(K)= 2g(K)+b(K)+|K|-2$.
\end{defn}

Next, see that the weaker of the lower bounds is additive under composition:

\begin{prop}
Let K be a composite link with the factor links $K_{1},...,K_{n}$.
Then \[
c(K)\geq \sum _{i=1}^{n}[2g(K_{i})+b(K_{i})+|K_{i}|-2]\]
\end{prop}
\begin{proof}
Use the facts that the genus is additive under composition (see e.
g. \cite{Ad}), and for the braid index it holds 
$b(K_{1}\#K_{2})=b(K_{1})+b(K_{2})-1$\cite{BirBraid}.
\end{proof}
\begin{rem*}
The free genus is also additive under knot composition, see \cite{freegenus}.
Thus for knots, we may also write here $g_{f}(K)$ instead of $g(K)$,
sometimes leading to better results.
\end{rem*}

\begin{cor}
If K is a composite link with factor links $K_{1},...,K_{n}$, and
all factor links are in $\mathcal{F}$, then
\[c(K)=\sum _{i=1}^{n}c(K_{i})\]
\end{cor}

At this stage, we already see that the crossing number is additive for torus
links; a problem repeatedly posed by Adams (see \cite{Ad}, chapter 5.1); 
we see that the answer
is almost as simple as the question in this special case.

With two other results by Murasugi, we have that $\mathcal{F}$ also
contains other classes of links:
We repeat here almost exactly corollary 2 in \cite{MuBraid}:

\begin{prop}
\cite{MuBraid} If K is an alternating fibered link, then K is in $\mathcal{F}$.
\end{prop}
\begin{proof}
Murasugi's statement is $c(K)=deg\Delta _{K}+b(K)-1$, and with $deg\Delta _{K}-|K|+1\leq 2g(K)$,
the rest follows.
\end{proof}

The following proposition is very closely related to proposition 7.4
in \cite{MuBraid}. However, I redo it here to point out a different
aspect: to fit the definition of the family $\mathcal{F}$.

\begin{prop}
Let $\gamma $ be a homogeneous n-braid, and L be the closure of $\gamma $.
If $b(L)=n$, then L is in $\mathcal{F}$.

\end{prop}
\begin{proof}
Draw a diagram D by simply closing $\gamma $. You have for the degree
of the reduced Alexander polynomial $deg\Delta _{L}=c(D)-s(D)+1$,
since $\gamma $ is homogeneous. With $deg\Delta _{L}+1-|K|\leq 2g(L)$
, $s(D)=b(L)$, this leads to $2g(L)\geq c(D)-b(L)-|K|+2$, and because
of $g(D)\geq g(L)$, $c(D)=2g(L)+b(L)+|K|-2$. By Lemma \ref{bound1}, you see
that c(D) is minimal.
\end{proof}

\section{Examples.}
In the previous chapter, we have seen some families of $\mathcal{F}$-links.
You may have observed that all links presented there were
homogeneous (for a definition see\cite{Cro}), and fibered.
However, it is neither true that every fibered homogeneous link is an
$\mathcal{F}$-link nor that every $\mathcal{F}$-link is fibered.
We take a look at an example:
\begin{xmpl}
The Perko knot (denoted by $10_{161}$ in Rolfsen's table \cite{Ro})
is fibered and homogeneous, and its (canonical) genus is 3, as well as
its braid index. But its crossing number is 10. Therefore, the Perko knot is
not in $\mathcal{F}$.
\end{xmpl}

I proceed with a discussion of alternating $\mathcal{F}$-links.
\begin{defn}
A number $n(D)$ associated with an alternating link diagram D is called an
alternating link invariant if $n(D_1)=n(D_2)$ for any two reduced alternating
diagrams $D_1,D_2$ of the same link.
\end{defn}
Consider again Seifert circles:
\begin{prop}
Let $D(L)$ be an alternating link diagram. Let $s_a(D)$
be the number of Seifert circles in $D$. Then $s_a(D)=:s_a(L)$ is an
alternating link invariant.
\end{prop}
\begin{proof}
Let $D$ be a reduced alternating link diagram of L. Then every
other reduced alternating diagram of L can be reached from $D$ by a
finite number of flyping moves \cite{MuFlype}.
The figure shows the flyping move.
\center{\includegraphics[width=1.0\textwidth]{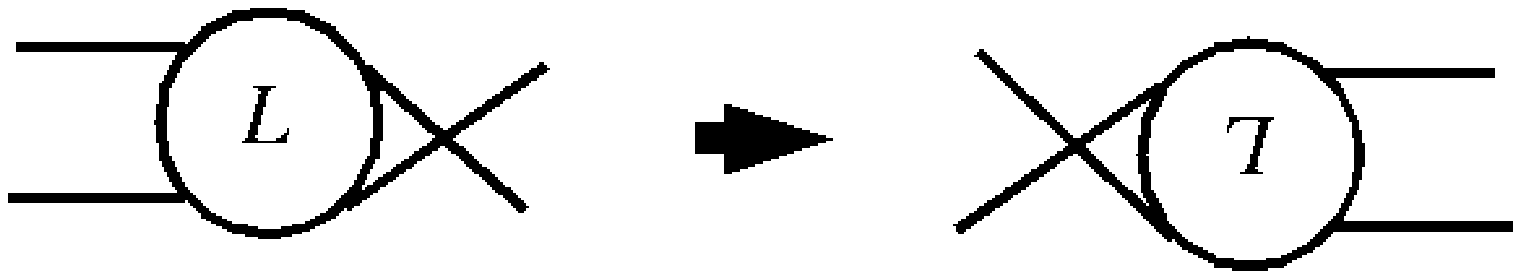}}
Hence, it is enough to show that $s(D)$ is not altered by the flyping move.
When considering only the underlying Seifert circles in the plane,
then flyping does not affect the connectedness
of the circles inside the disk. If the strands on the right-hand side
of the above figure have the same orientation, the flyping move
does not change the number of Seifert circles.
For the remaining case, we use the idemposition
algebra for arrangements of
(unoriented) circles in the plane (see\cite{GSB}). In this case, the effect
of flyping outside the disk looks like:
\center{\includegraphics[width=1.0\textwidth]{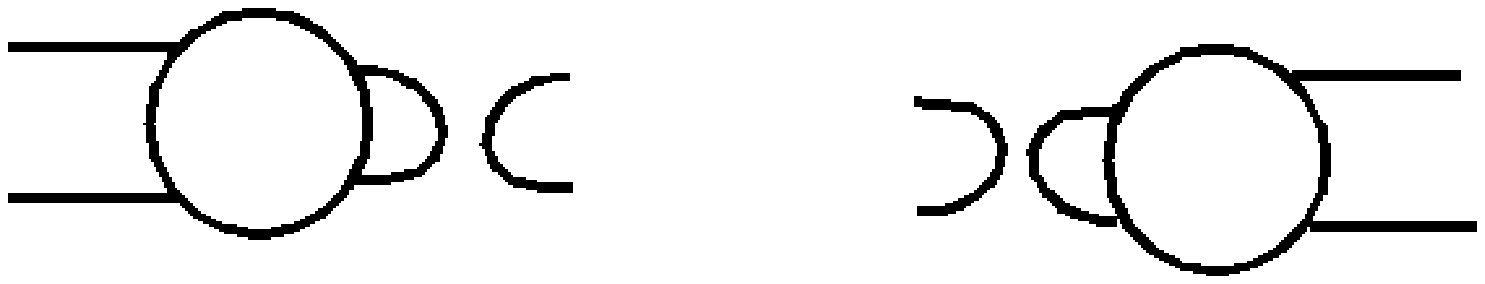}}

When counting Seifert circles, we must consider several cases
for the above disk diagram on the left:
The outgoing arcs of the diagram can be connected outside the diagram
in two ways, and there are two possibilities how the ingoing arcs
into the disk can be connected inside the disk:
\begin{enumerate}
\item the upper left arc is connected to the lower
left arc outside the diagram (called numerator closure)
	\begin{enumerate}
	\item the upper left arc is connected inside the disk to the
	lower left arc
	\item the upper left arc is connected inside the disk to the right ingoing
	arc
	\end{enumerate}
\item the upper left arc is connected to the upper right arc
outside the diagram (called denominator closure)
	\begin{enumerate}
	\item the upper left arc is connected inside the disk to the
	lower left arc inside the disk
	\item the upper left arc is connected to the right ingoing arc inside the disk.
	\end{enumerate}
\end{enumerate}
Now the important thing is to see that the left arcs are in all four cases
connected to each other.
We further need two idemposition types, see the following picture:
\center{\includegraphics[width=1.0\textwidth]{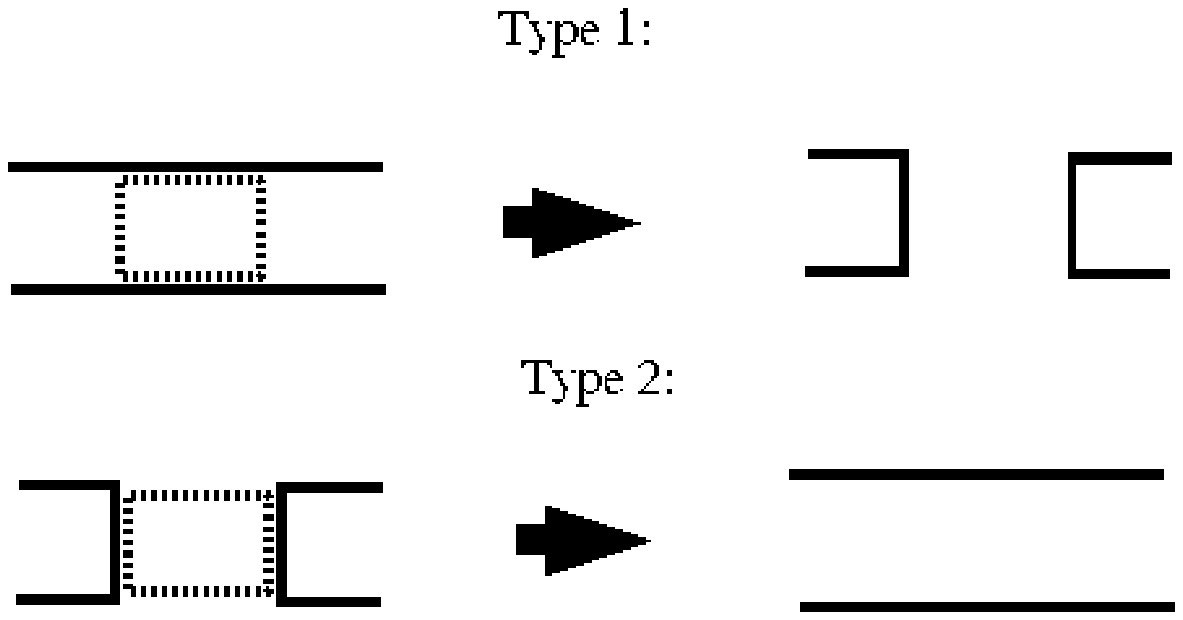}}
To imitate flyping, do the following:
Apply first a type 1 idemposition $\Phi$ on the left side.
In all cases, theorem 14 in \cite{GSB} gives $s(\Phi(D))=s(D)+1$.
The disc is now disconnected from the rest of the link.
Turn it in the plane by $\pi$.
Next, apply a type 2 idemposition $\Theta$
on the right side. Since the arcs on the right side of the image
cannot be connected before applying $\Theta$, again by theorem 14,
we have $s(\Theta(\Phi(D)))=(s(D)+1)-1$, and the proof is complete.
\end{proof}
We obtain immediately the following condition for alternating
$\mathcal{F}$-links:
\begin{cor}\label{sab}
An alternating link is in $\mathcal{F}$ iff $s_a(L)=b(L)$.
\end{cor}
\begin{proof}
Assume first $s_a(L)=b(L)$.
It is a long known fact that the Seifert surface from an alternating diagram
$D$ is minimal, i. e. $g(D)=g(L)$.  Since $s_a(L)=s(D)=b(L)$, and a reduced
alternating diagram has minimal number of crossings \cite{Ad},
$c(L)=2g(L)+b(L)-1$.
On the other hand, if $b(L)<s_a(L)$, for a reduced alternating diagram we have
$c(D)>2g(D)+b(L)-1$.
\end{proof}
We proceed with a statement on planar Murasugi sums. For the definition
of the terms used here, see \cite{MuBraid}.
\begin{cor}
Let L be an alternating link, and $D(L)$ be the planar Murasugi sum
of special alternating
link diagrams $D(L_1),..., D(L_n)$.
then
\[s_a(L)-1 = \sum_{i=1}^n (s_a(L_i)-1)\]
\end{cor}
\begin{proof}
Follows immmediately by the definition of planar Murasugi sum.
\end{proof}
It is now natural to ask whether a similar equality holds
for the braid index. This is conjectured in \cite{MuBraid}:
\begin{conj}\label{MurasugiSum}
Let $L$ be the planar Murasugi sum of special alternating
links $L_1,..., L_n$.
then
\[b(L)-1 = \sum_{i=1}^n (b(L_i)-1)\]
\end{conj}
I note that the special alternating $\mathcal{F}$-links are exactly those
satisfying theorem 8.1 in \cite{MuBraid}, called there
"special alternating links of the nonmultiple type". This is an 
important point for the following discussion.

\begin{xmpl}
The pretzel links $\mathcal{P}(a_1,a_2,...a_{2n})$ with
all $a_i\geq 2$ are special $\mathcal{F}$-links
(which are not fibered).
\end{xmpl}

\begin{prop}
Let $D(L)$ be the planar Murasugi sum of special alternating
link diagrams $D(L_1),..., D(L_n)$. Assume at least one of
$L_1,..,L_n$ is a non-$\mathcal{F}$-link.
Then $L$ is a non-$\mathcal{F}$-link.
\end{prop}
\begin{proof}
In the proof of theorem 8.1 in \cite{MuBraid}, a local
diagrammatic move upon a special alternating diagram of a non-$\mathcal{F}$
link reduces the number of Seifert circles.
Use this move in $D(L)$ to obtain a diagram $D'$ with fewer
Seifert circles. But now, $b(L)\leq s(D') < s_a(L)$. Recall \ref{sab}
to conclude that L is not in $\mathcal{F}$.
\end{proof}
\begin{cor}
Let $L$ be a $\mathcal{F}$-link.
Then conjecture \ref{MurasugiSum} holds.
\end{cor}
\begin{proof}
Since $L$ is in $\mathcal{F}$, $b(L)=s_a(L)=\sum_{i=1}^n (s_a(L_i)-1)$.
On the other hand, for all $L_i$, $s_a(L_i)=b(L_i)$ because
they must also be in $\mathcal{F}$ (which follows from the above proposition).
\end{proof}

I believe furthermore the alternating links constructed from
special alternating $\mathcal{F}$-links are exactly
the alternating $\mathcal{F}$-links.
This would follow immediately from conjecture \ref{MurasugiSum}
and complete the classification of alternating
$\mathcal{F}$-links.

\section{Application to satellite knots.}

We can apply Lemma \ref{bound1} to give a "good" estimate for
the crossing number of (p,q)-cable knots about a knot in $\mathcal{F}$.
(Here, p means the linking number of the knot with a meridian of 
the essential knotted torus in which L lies.)

\begin{cor}
Let p,q be positive integers with $\gcd (p,q)=1$.
Let K be the (p,q)-cable knot about a $\mathcal{F}$-knot C.
Then $c(K)\geq q(p-1) + p \cdot c(C)$.
\end{cor}
\begin{proof}
Schubert showed in (\cite{Schubert}, pp. 247 seq.) that
\[ 2 g(K)= (p-1)(q-1)+2p \cdot g(C) \]
Furthermore, the braid index of K is $p \cdot b(C)$ \cite{BirSat, WillBraid}.
Since $\gcd (p,q)=1$, $|K|=1$. Use \ref{bound1} to obtain
$c(K)\geq q(p-1) + p \cdot c(C)$.
\end{proof}
\begin{rem*}
However, this bound is far from being
sharp, as we conjecture that the crossing number is
$q(p-1) + p^2 \cdot c(C)$. This may illustrate the difficulties
in proving statements about the crossing number of satellite knots.
\end{rem*}
Similar estimates can be easily established in the same spirit,
which cover other special cases, as you may figure out yourself.
Nevertheless, we state here another example.

(The definitions for pattern types and the weights $w_i$
used in the following are the same as in \cite{Nutts,BirSat}.)

\begin{cor}
Let K be a satellite knot with a $\mathcal{F}$-knot companion C
and a closed alternating braid B as type 0 pattern.
Then, $c(K) \geq c(B)+ b(B)\cdot c(C)$.
\end{cor}
\begin{proof}
Another result of Schubert (\cite{Schubert},p.192) states that in our case
\[2g(K) \geq 2b(B) \cdot g(C) + 2g(B)\]
Since B is an $\mathcal{F}$-knot, $2g(B) = c(B)-b(B)+1$; the same holds for C.
Here again $b(K)= b(B)\cdot b(C)$ holds, for B is a type 0
pattern \cite{BirSat}.
Now, \ref{bound1} together with both B and C being $\mathcal{F}$-knots:
\[c(K) \geq c(B)+ b(B)\cdot c(C)\]
\end{proof}

We can also show an unconditional estimate for satellite knots, as follows:

\begin{prop}\label{sat1}
Let K be a satellite knot with companion C. Let B be the
knot from the pattern. Then
\[ c(K) \geq
\left\{
\begin{array}{l}
w_0 \cdot(2g(C) + b(C) -  1) + g(B) + w_0 \textnormal{ for pattern type 0} \\
w_0 \cdot(2g(C) + b(C) -  1) + g(B) + w_0 + w_1 \textnormal{ for pattern type 1}\\
2\alpha(C)-2 \textnormal{ else}\\
\end{array}
\right.\]
There, $\alpha(K)$ is the arc index of K as defined in \cite{CroArc}.
\end{prop}
\begin{proof}
First consider type 0 and type 1 patterns:
Schubert's theorem on the genus of satellite knots
reads for the general case as follows:
$g(K) \geq n \cdot g(C) + g(B)$, where B is the knot from the pattern
n is the linking number of a meridian of the essential torus with B.
Here, $n=w_0$.
The braid index of K is $w_0 \cdot b(C)$ for a type 0 pattern
and $w_0 \cdot b(C)+w_1$ for a type 1 pattern respectively \cite{BirSat}.
Using Lemma \ref{bound1}, we get
$c(K)\geq w_0 \cdot(2g(C) + b(C) -  1) + g(B) + w_0$ for a type 0 pattern,
resp. $c(K)\geq w_0 \cdot(2g(C) + b(C) -  1) + g(B) + w_0 + w_1$
for a type 1 pattern.

For the other case, we use a different estimate given by Ohyama \cite{Ohyama}:
$c(K)\geq 2b(K)-2$. Since we have a type k pattern,
$b(K)\geq \alpha(C)$ (\cite{Nutts}, prop.3.3.5).
We can conclude $c(K)\geq 2 \alpha(C)-2$ in this case.
\end{proof}

We can furthermore solve here the problem whether the crossing number
of a satellite link is greater than that of the companion
(\cite{Kirby}, problem 1.67) if the latter
is an alternating fibered knot:
\begin{cor}
Let K be a satellite knot and let its companion C be an alternating fibered
knot. Then $c(K)\geq 2c(C)+2$.
\end{cor}
\begin{proof}
Cromwell showed that $\alpha(K)\geq c(K)+2$ if K is an alternating knot
\cite{CroArc}.
Since $K \in \mathcal{F}$, and for a proper satellite knot,
the weight $w_0$ is greater than 1, the rest follows from the above
proposition.
\end{proof}
The result can be generalized to more companion types and sharpened
at the same time if we restrict the pattern type:
\begin{cor}
Let K be a satellite knot and let its companion C be a $\mathcal{F}$-knot.
Furthermore, assume K has a pattern of type 0 or 1.
Then $c(K)> w_0 \cdot c(C)$.
\end{cor}
\begin{proof}
Obvious from \ref{sat1}.
\end{proof}

The results worked out here are all dealing with the Seifert genus.
But the canonical genus gives us more power;
we deal with this in the next section. However, I don't believe
that any satellite link is in $\mathcal{F}$.

\section{Crossing number, canonical genus and knot polynomials}

Of course, the question of determining the minimal crossing number
remains open in general. But we can state another
beautiful lower bound in terms of the HOMFLY polynomial, which is
again additive under composition. Define the two-variable Laurent
polynomial in the variables $v,z$ $P[v,z]=P(K)$ as in \cite{HOMFLY},
and write shortly $e$ (resp. $m$) and $E$ (resp. $M$) for
the minimal and maximal non-zero exponent in $v$ (resp. $z$).

\begin{cor}
For every link $K$,
\[c(K)\geq M+\frac{1}{2}(E-e)\]
\end{cor}
\begin{proof}
Use Lemma \ref{bound1} together with $M\leq \min _{D=D(K)}{c(D)-s(D)+1}$
\cite{MoSeif} and the MWF- inequality
\cite{MoSeif,FranksWilliams}:
$b(K)\geq \frac{1}{2}(E-e)+1$
Furthermore, $|K|=m$.
\end{proof}
Notice that the right-hand side is again additive under composition,
since for a composite link $L\#M$ the polynomial is $P(L)\cdot P(M)$.

\begin{rem*}
Sometimes, this bound can lead to better results than
the inequality using the genus. For instance, $M=6$ for the untwisted
double of the trefoil, whereas the free genus is 2, and the genus
is only 1.
On the other hand, Stoimenow showed an example where $2g>M$
\cite{Sto}; and there are also some examples for which
the MWF- inequality is not sharp.
\end{rem*}
We continue with showing up a connection between the crossing number
of a knot and the canonical genus of its double.
(I noticed this when reading in \cite{MoSeif}
that the canonical genus of the Whitehead double of the trefoil is 3,
which is also the minimal crossing number)

\begin{prop}\label{bound2}
Let K be a knot and $W_K$ be a Whitehead double of K. Then
\[c(K)\geq \widetilde{g}(W_K)\]
\end{prop}

\begin{proof}
First, assume that $W_K$ is untwisted.
The proof is constructive: Draw a diagram with minimal crossing number
of K. Then choose a parallel to run along the curve in the inverse
direction of the original curve
(diagrams of this kind are often called blackboard diagrams).
Then, in a region which is not near
a crossing, replace the antiparallel by a (say) positive clasp
to obtain a diagram D of $W_K$ with $4c(K)+2$ crossings.

Cut out the crossings to obtain a bunch of Seifert circles with a
pattern as shown in the figure.

\includegraphics[width=1.0\textwidth]{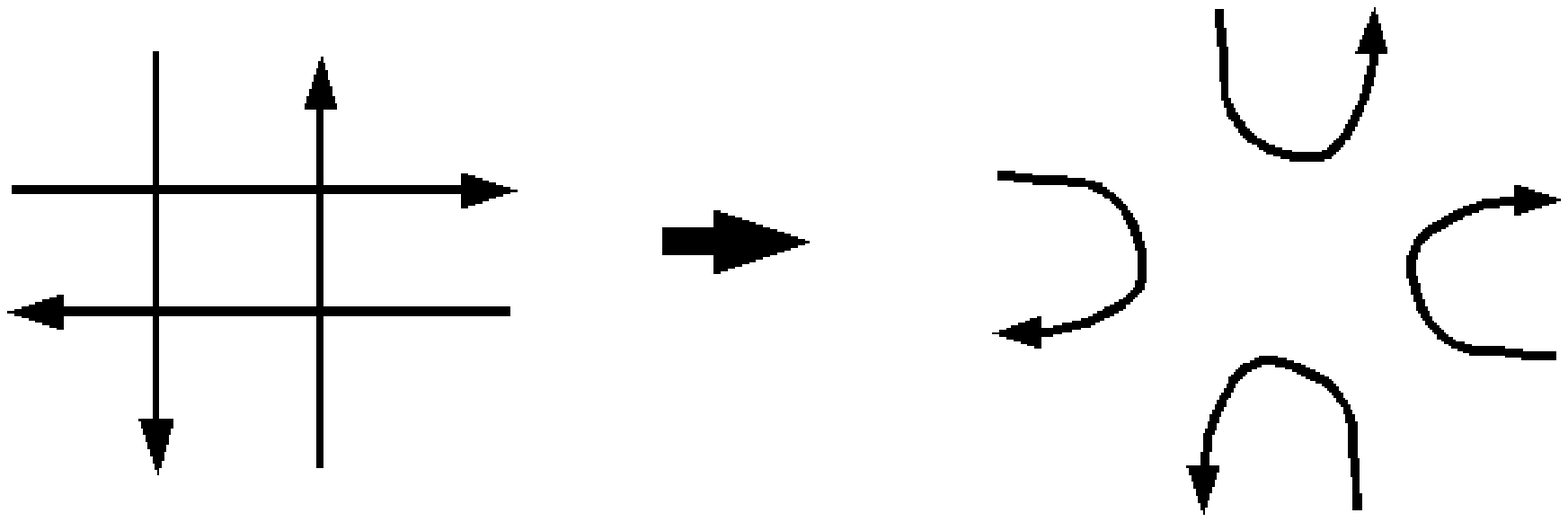}

Near each underlying crossing of K, you have four emerging arcs. Since
the underlying diagram of K is minimal (and hence reduced),
every one of these four arcs belongs to
a distinct Seifert circle. In the region with the clasp from
the Whitehead doubling, you see a closed Seifert circle
and two lines above and below it, belonging to the same Seifert circle.
The following picture should illustrate this.

\includegraphics[width=1.0\textwidth]{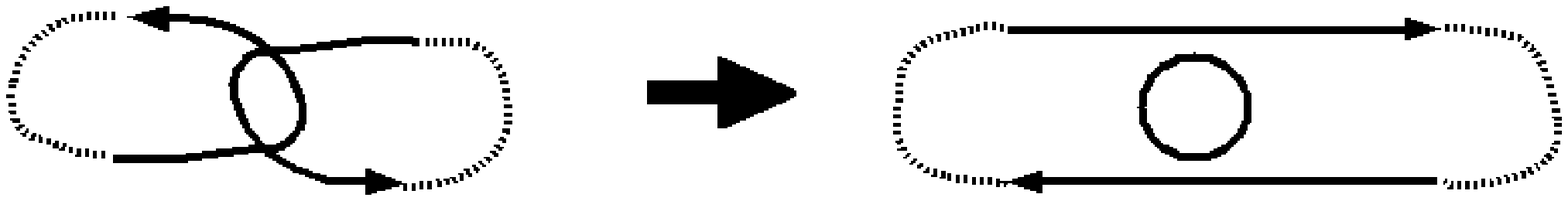}

Since every arc closes up to a Seifert circle with another arc,
we count $s(D)=2c(K)+1$ Seifert circles in the diagram, and
clearly, this leads for the genus of this diagram 
to $g(D)=c(K)\geq \widetilde{g}(W_K)$ .

Next, consider the case $W_K$ is twisted.
Let $D_n$ (resp. $D_{n+1}$) be a diagram with n (resp. $n+1$)
half-twists.
The following picture illustrates that one more half-twist
produces splits a Seifert circle in $D_n$ into two, and at the same
time we have more crossing in $D_{n+1}$, thus every
pair of diagrams $D_n$, $D_{n+1}$
has the same genus. See the figure.

\includegraphics[width=1.0\textwidth]{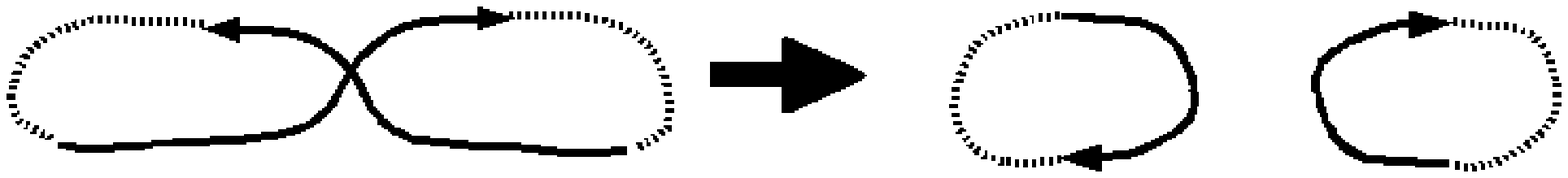}

Complete the proof by induction on n.
\end{proof}

\begin{rem*}
I was informed that this proposition was found 
independently by Kidwell and Stoimenow \cite{StoKid}.
\end{rem*}

Together with Morton's canonical genus inequality, we can
estimate the minimal crossing number from the polynomial
of a satellite about it. A similar technique has already successfully
been used in \cite{MoCable} for the braid index in special cases.

\begin{problem}
Does an infinite family of knots exist with $\widetilde{g}(W_K)=c(K)$?
\end{problem}

Very recently, Stoimenow and Kidwell asked a question
(see \cite{suki}, Problem 1.15) which turns
out to be closely related to this problem:
\begin{conj}
Let K be a knot and $W_K$ be its Whitehead double. Let
$F(K)$ be the Kauffman polynomial of K in the variables $a$ and $z$.
Then
\[ 2(\max\deg_z F(K) +1)=\max\deg_z P(W_K)\]
\end{conj}

The truth of the conjecture would imply that there exists
such an infinite family, namely the prime alternating knots.
With the inequality found here, I can give
further "evidence"
for the truth of the above equality
(both sides are bounded above by $2c(K)$), and
attack a part of it:

\begin{prop}
Let K be a prime alternating knot and Wh(K) a
Whitehead double of K. Then
\[ 2(\max\deg_z F(K) +1)\geq \max\deg_z P(W_K)\]
\end{prop}
\begin{proof}
Kidwell showed that $c(K) \geq \max\deg_z F(K) +1$ with equality if K is a
prime alternating knot \cite{Kidwell}. Morton's inequality and
proposition \ref{bound2} give us
$2c(K)\geq \max\deg_z P(W_K)$, and the proof is complete.
\end{proof}

Since we have also found now
\[2c(K)\geq \max\deg_z P(W_K)\]
and
\[2c(K)\geq 2(\max\deg_z F(K) +1)\]
I feel free to mention the following, which is quite similar:

Let $G(K)$ be the Kauffman polynomial with coefficients reduced modulo 2.
Then $b(W_K)\geq \alpha(K)$ (see \cite{Nutts}, prop.3.3.5) together with
$\alpha(K) \geq \max\deg_a G(K) - \min\deg_a G(K)+2$
(see \cite{Nutts}, prop.4.4.1) imply:
\[ 2b(W_K)-2\geq 2(\max\deg_a G(K) - \min\deg_a G(K)+1)\]
whereas
\[2b(W_K)-2 \geq \max\deg_v P(W_K) -\min\deg_v P(W_K)\]
is the MWF-inequality. (Note that the degrees in $v$ of the HOMFLY polynomial
depend heavily on how $W_K$ is twisted, whereas
$G(K)$ is of course not affected.)

\begin{acknowledgement*}
I would like to thank Joan Birman, and Alexander Stoimenow for answering my
questions, and providing useful hints.
\end{acknowledgement*}

\end{document}